\documentclass[11pt, a4paper]{amsart}

\usepackage{amsmath, amsthm, amssymb}
\usepackage{txfonts, mathrsfs}
\usepackage[usenames]{color}
\usepackage{psfrag}
\usepackage{graphicx}
\usepackage{todonotes}
\usepackage{comment}
\usepackage[normalem]{ulem}
\usepackage[hidelinks, pdfencoding=auto, psdextra]{hyperref}
\usepackage{cleveref}
\usepackage{microtype}

\usepackage[final]{showlabels}

\numberwithin{equation}{section}

\newtheorem{thm}{Theorem}[section]
\newcommand{\bt}{\begin{thm}}
\newcommand{\et}{\end{thm}}

\newtheorem{cor}[thm]{Corollary}
\newcommand{\bc}{\begin{cor}}
\newcommand{\ec}{\end{cor}}

\newtheorem{lem}[thm]{Lemma}
\newcommand{\bl}{\begin{lem}}
\newcommand{\el}{\end{lem}}

\newtheorem{prop}[thm]{Proposition}
\newcommand{\bp}{\begin{prop}}
\newcommand{\ep}{\end{prop}}

\newtheorem{defn}[thm]{Definition}
\newcommand{\bd}{\begin{defn}}
\newcommand{\ed}{\end{defn}}

\newtheorem{rmrk}[thm]{Remark}
\newcommand{\br}{\begin{rmrk}}
\newcommand{\er}{\end{rmrk}}

\newtheorem{quest}[thm]{Question}
\newcommand{\bq}{\begin{quest}}
\newcommand{\eq}{\end{quest}}

\newtheorem{example}[thm]{Example}

\theoremstyle{definition}

\newcommand{\R}{\mathbb{R}}

\newdimen\vintkern\vintkern12pt
\def\vint{-\kern-\vintkern\int}

\newcommand{\hm}{{\mathcal H}}

\newcommand{\trace}{\operatorname{tr}}
\newcommand{\length}{\ell}

\newcommand{\md}{\operatorname{md}}

\newcommand{\ap}{\operatorname{ap}}
\newcommand{\apmd}{\ap\md}
\DeclareMathOperator{\MOD}{mod}

\begin{document}

\title{Monotone Sobolev extensions in metric surfaces and applications to uniformization}

\author{Damaris Meier}
\address
  {Department of Mathematics\\ ETH Zurich \\ R\"amistrasse 101\\ 8092 Zurich, Switzerland}
\email{damaris.meier@math.ethz.ch}

\author{Noa Vikman}
\address
  {Department of Mathematics\\ University of Fribourg\\ Chemin du Mus\'ee 23\\ 1700 Fribourg, Switzerland}
\email{noa.vikman@unifr.ch}

\author{Stefan Wenger}
\address
  {Department of Mathematics\\ University of Fribourg\\ Chemin du Mus\'ee 23\\ 1700 Fribourg, Switzerland}
\email{stefan.wenger@unifr.ch}

\date{\today}

\thanks{Research supported by Swiss National Science Foundation Grant 212867}

\begin{abstract}
We prove a monotone Sobolev extension theorem for maps to Jordan domains with rectifiable boundary in metric surfaces of locally finite Hausdorff 2-measure. This is then used to prove a uniformization result for compact metric surfaces by minimizing energy in the class of monotone Sobolev maps.
\end{abstract}

\maketitle

\section{Introduction}

Every rectifiable Jordan curve $\Gamma$ in the Euclidean plane can easily be parametrized by a Lipschitz curve $\gamma\colon S^1\to\Gamma$, for example through the constant speed parametrization. It is well known that $\gamma$ then extends to a Lipschitz map from the closed unit disc $\overline{D}$ to the Jordan domain bounded by $\Gamma$. If the Lipschitz map $\gamma$ is moreover monotone, meaning that it is continuous and the preimage of every point is connected, then there exists a Lipschitz extension of $\gamma$ that is also monotone. The shortest curve extension, e.g.\ as defined in \cite[Section 6]{Hencl-Koski-Onninen_3-dim}, can serve as such an extension. In this note we are interested in generalizations of monotone Lipschitz extension results to a non-smooth setting. 

\subsection{Main results} We consider the case where the target space $X$ is a metric surface, i.e.\ a complete metric space homeomorphic to a 2-dimensional surface, possibly with non-empty boundary. Singular metric surfaces naturally arise in various mathematical contexts and applications; for example as boundaries, limits or deformations of classical smooth objects. We assume that $X$ has locally finite Hausdorff 2-measure, for a definition we refer to \Cref{sec:notions}. Note that such a space $X$ can contain a purely 2-unrectifiable part that is dense in $X$. In particular, there might not exist any Lipschitz maps from $\overline{D}$ to $X$ covering an open set in $X$. Thus, there is no hope for extending any map from $S^1$ to $\partial\Omega$ to a monotone Lipschitz map from $\overline{D}$ to $\overline{\Omega}$, where $\Omega$ is a Jordan domain in $X$. Instead, we look for extensions within the class of monotone Sobolev mappings. Note that Sobolev maps from Riemannian (or more general) domains into complete metric spaces have been an important object of study within the past 35 years, and by now a robust theory has emerged, see \cite{KS93,Jost94,Haj96,HKST15}. We state the following main theorem of this note.

\bt \label{thm:monotone-filling}
    Let $X$ be a metric surface with locally finite Hausdorff 2-measure. Let $\Omega\subset X$ be a Jordan domain with rectifiable boundary. Then every monotone map in $W^{1,2}(S^1,\partial \Omega)$ extends to a monotone map in $W^{1,2}(\overline D,\overline \Omega)$.
\et

For the definitions of Sobolev spaces, we refer to \Cref{sec:Sobolev}. Note that in the above theorem we endow $\partial \Omega$ and $\overline{\Omega}$ with the metric inherited from $X$. We will prove \Cref{thm:monotone-filling} in \Cref{sec:Sobolev-monotone-extension} first under the additional assumption that the metric surface $X$ is locally geodesic. The proof  utilizes the existence and regularity of energy minimizing Sobolev maps from \cite{MW21}, which builds on \cite{LW-area,lytchak_energy_2017,LW20}. Moreover, we make use of a collar construction similar to that in \cite{Stadler21} or \cite{Creutz22}. By exploiting the more recent work \cite{NR22} instead of \cite{MW21} one can drop the condition of $X$ being locally geodesic, see \Cref{rmk:NR22}. In a next step we establish a monotone Sobolev extension result, \Cref{thm:monotone-sobolev-extension-compact-surface}, for domains of higher topology by decomposing the surfaces into suitable Jordan domains and applying \Cref{thm:monotone-filling}.

A result of Youngs \cite{youngs_homeomorphic_1948} implies that surjective monotone mappings between compact homeomorphic 2-manifolds are exactly the uniform limits of homeomorphisms. Improving the statement of \Cref{thm:monotone-filling} to obtain homeomorphic Sobolev extensions is impossible for such general targets as in the theorem, even if the map from $S^1$ to $\partial\Omega$ is a Lipschitz embedding, see \Cref{ex:collapsed-disc}. This stands in contrast to recent developments for homeomorphic Sobolev extensions in the Euclidean plane by Koski, Onninen and co-authors \cite{Koski-Onninen_JEMS,Koski-Onninen_Sobolev-Schoenflies,Koski-Onninen_Bi-Sobolev,Hencl-Koski-Onninen_3-dim,Koski-Onninen-Xu_Internal}.

We now provide an application of \Cref{thm:monotone-filling} to the uniformization of metric surfaces. The classical uniformization theorem states that every simply connected Riemann surface is conformally equivalent to the unit disc, the complex plane or the Riemann sphere. Since the beginning of this century, the problem of finding analogues of this statement for non-smooth metric surfaces has drawn a lot of research interest, especially due to its relevance to geometric group theory and complex dynamics, see e.g.\ \cite{BK02,BK05,Raj17,LW20,MW21,NR21,NR22,Mei24}. A natural approach to producing a uniformization map is by minimizing energy in a certain class of Sobolev mappings spanning the given metric surface. Following this strategy, we will prove the next result.

\bt \label{thm:uniformization}
    Let $X$ be a locally geodesic metric space that is homeomorphic to a compact orientable smooth surface $M$, has rectifiable boundary (if non-empty), and is of finite Hausdorff 2-measure. Then the family $\Lambda$ of surjective and monotone Sobolev maps from $M$ to $X$ is non-empty. Further, there exists $u\in \Lambda$ and a Riemannian metric $g$ on $M$ such that
    $$
        E_+^2(u,g)=\inf\{E_+^2(v,h): v\in \Lambda, \ h \text{ Riemannian metric on } M\}.
    $$
    For any such pair $(u,g)$ the map $u\colon (M,g)\to X$ is infinitesimally isotropic, and, in particular, yields a weakly $K$-quasiconformal parametrization of $X$ with the optimal constant $K=4/\pi$. 
\et

Here, $E^2_+(u,g)$ denotes the Reshetnyak energy of $u$ with respect to $g$, see \Cref{sec:Sobolev}. Infinitesimal isotropy is a strong metric variant of weak conformality. A definition of infinitesimal isotropy can be found in \Cref{sec:metric-diff}. 
A weakly $K$-quasiconformal parametrization of $X$ is a monotone and surjective map $u\colon M\to X$ satisfying the distortion inequality 
\begin{equation*}
        \MOD(\Gamma)\leq K\cdot \MOD(u\circ\Gamma)
\end{equation*}
for every family $\Gamma$ of curves in $M$, where $\MOD$ refers to the conformal modulus of curve families as defined in \Cref{sec:modulus}. Note that under the additional assumption of $X$ being reciprocal, every weakly quasiconformal parametrization of $X$ upgrades to being a geometrically quasiconformal homeomorphism, see \cite[Section 3]{MW21} and \cite[Section 7]{NR22}. 

The Riemannian metric $g$ in the theorem can be chosen to be of constant sectional curvature $-1, 0$ or $1$ and such that $\partial M$ is geodesic (if non-empty). And a smooth surface refers to a smooth, connected $2$-dimensional manifold, possibly with non-empty boundary. The first statement of \Cref{thm:uniformization} follows as a corollary of the monotone Sobolev extension results \Cref{thm:monotone-filling} and \Cref{thm:monotone-sobolev-extension-compact-surface}. To show the existence of an energy minimizing pair $(u,g)$ we follow a direct variational method.

\Cref{thm:uniformization} provides a weakly quasiconformal parametrization result for any locally geodesic compact metric surface with possibly empty rectifiable boundary and finite Hausdorff 2-measure by following the natural approach of energy minimization. This was previously achieved for metric discs in \cite{MW21} and metric surfaces with non-empty boundary in \cite{Mei24}. Hence, the novelty of \Cref{thm:uniformization} is the parametrization of closed surfaces via energy minimizers. 

By following a completely different proof strategy, Ntalampekos and Romney \cite{NR22}, see also \cite{NR21}, were able to produce a weakly quasiconformal uniformization result, where the only assumption on the metric surface $X$ is the local finiteness of the Hausdorff 2-measure. Weakly quasiconformal uniformization already enjoys a great variety of applications in fields such as geometric function theory, geometric measure theory and complex analysis, see e.g.\ \cite{EIR23,MN24,BMW25,MR_finite-distortion,MR-quasiconformality,Nta_conf-remove,Nta_smooth-approx}.

\subsection{Open questions}
Recent work of Koski, Onninen and Xu \cite{Koski-Onninen-Xu_Internal} shows that if $\Omega$ is a Jordan domain in the Euclidean plane, then a homeomorphism $\gamma\colon S^1\to\partial\Omega$ extends to a homeomorphism from $\overline{D}$ to $\overline{\Omega}$ in the Sobolev class $W^{1,2}(D,\overline\Omega)$ if and only if $\gamma$ satisfies the so-called internal 2-Douglas condition
\begin{equation}\label{eq:intrinsic-Douglas}
    \int_{S^1}\int_{S^1}\frac{d_i(\gamma(x),\gamma(y))^2}{|x-y|^2}\,dx\,dy<\infty,
\end{equation}
where $d_i$ denotes the intrinsic metric in $\overline{\Omega}$. In other words, $\gamma$ is in the fractional Sobolev space $W^{\frac{1}{2},2}(S^1,(\overline\Omega,d_i))$, see e.g.\ \cite{Chiron07}. We do not know whether this weaker boundary regularity is sufficient for the existence of a monotone Sobolev extension in a metric surface setting.
\bq
    Does the statement of \Cref{thm:monotone-filling} still hold if we impose that the monotone boundary map satisfies \eqref{eq:intrinsic-Douglas} instead of being in $W^{1,2}(S^1,X)$?
\eq
Recall that within the setting of \Cref{thm:monotone-filling}, we can not hope for a homeomorphic Sobolev extension result in the sense that every Sobolev homeomorphism from $S^1$ to $\partial \Omega$ extends to a Sobolev homeomorphism from $\overline D$ to $\overline{\Omega}$, see \Cref{ex:collapsed-disc}. The above mentioned fact that the map $u$ from \Cref{thm:uniformization} upgrades to being a homeomorphism after 
assuming that $X$ is reciprocal, as well as the results \cite{Koski-Onninen_JEMS,Koski-Onninen_Sobolev-Schoenflies,Koski-Onninen_Bi-Sobolev,Hencl-Koski-Onninen_3-dim,Koski-Onninen-Xu_Internal}  on homeomorphic Sobolev extensions in the Euclidean plane motivate the next question.
\bq
    What are the minimal assumptions on a metric surface under which a homeomorphic Sobolev extension result holds?
\eq

The structure of the article is as follows: In \Cref{sec:prelim} we introduce notation and give some background on metric space valued Sobolev maps. \Cref{sec:Sobolev-monotone-extension} is devoted to proving \Cref{thm:monotone-filling} and its generalization, \Cref{thm:monotone-sobolev-extension-compact-surface}, to domains of higher topology, whereas \Cref{thm:uniformization} is proven in \Cref{sec:uniformization}.\\

\textbf{Acknowledgments:} We thank Dimitrios Ntalampekos, Kai Rajala and Matthew Romney for helpful comments on an earlier version of this note.

\section{Preliminaries}\label{sec:prelim}
\subsection{Basic notions}\label{sec:notions}
Let $(X,d)$ be a metric space. The open ball in $X$ centered at $x\in X$ of radius $r>0$ will be denoted by $B(x,r)$. The length of a curve $\gamma$ in $X$ is denoted by $\ell(\gamma)$. A curve $\gamma\colon[a,b]\to X$ is called rectifiable if $\ell(\gamma)<\infty$ and geodesic if $\length(\gamma) = d(\gamma(a),\gamma(b))$. The space $X$ is locally geodesic if every point has a neighborhood in which any two points can be joined by a geodesic with image in $X$.

We denote the $n$-dimensional Hausdorff measure of a set $\Omega\subset X$ by $\hm^n(\Omega)$. The normalizing constant is chosen in such a way that if $X$ is Euclidean $\R^n$, then $\hm^n$ agrees with the Lebesgue measure $\mathcal{L}^n$. If $M$ is a smooth $n$-dimensional manifold equipped with a Riemannian metric $g$, then the $n$-dimensional Hausdorff measure $\hm_g^n$ on $M$ coincides with the Riemannian volume. 

Recall from the introduction that a smooth surface is a connected, smooth $2$-dimensional manifold, possibly with non-empty boundary, and a complete metric space homeomorphic to a smooth surface is called a metric surface. Throughout this paper, the reference measure on a metric surface will always be the $2$-dimensional Hausdorff measure. 

Denote by $D$ and $S^1$ the open unit disc and unit circle in $\R^2$, respectively. Let $X$ be a metric surface. For a subset $\Omega\subset X$ we indicate the closure by $\overline \Omega$ and its topological boundary by $\partial \Omega$. A Jordan curve in $X$ is a subset of $X$ homeomorphic to $S^1$. Similarly, a Jordan arc in $X$ is a subset of $X$ homeomorphic to an interval. If $\Omega\subset X$  is homeomorphic to $D$ and $\partial\Omega$  is a Jordan curve, then we call $\Omega$ a Jordan domain in $X$.

\subsection{Conformal modulus and weak quasiconformality}\label{sec:modulus}
Let $X$ be a metric space and $\Gamma$ a family of curves in $X$. A Borel function $\rho\colon X\to [0,\infty]$ is said to be admissible for $\Gamma$ if $\int_\gamma \rho\geq 1$ for every locally rectifiable curve $\gamma\in\Gamma$. See \cite{HKST15} for the definition of the path integral $\int_\gamma\rho$. The (conformal) modulus of $\Gamma$ is defined by $$\MOD(\Gamma)= \inf_\rho\int_X\rho^2\,d\hm^2,$$ where the infimum is taken over all admissible functions for $\Gamma$. By definition, $\MOD(\Gamma)=\infty$ if $\Gamma$ contains a constant curve.

Assume that $X$ is homeomorphic to a compact smooth surface $M$. 
\bd\label{def:weak-qc}
    A monotone and surjective map $u\colon M\to X$ is called weakly $K$-quasiconformal parametrization of $X$ if 
    \begin{equation*}
        \MOD(\Gamma)\leq K\cdot \MOD(u\circ\Gamma)
    \end{equation*}
    holds for every family $\Gamma$ of curves in $M$, where $u\circ \Gamma$ denotes the family of all curves of the form $u\circ\gamma$ for some $\gamma\in\Gamma$.
\ed
Recall from the introduction that a map between topological spaces is monotone if it is continuous and the preimage of every point is connected. A result of Youngs \cite{youngs_homeomorphic_1948} implies that surjective monotone mappings between compact homeomorphic surfaces are exactly the uniform limits of homeomorphisms. For more properties of weakly quasiconformal mappings, we refer to \cite[Sections 2.4 and 7]{NR21}.

\subsection{Metric space valued Sobolev maps}\label{sec:Sobolev}
There exist several equivalent definitions of Sobolev maps from a Euclidean or Riemannian domain into a complete metric space, see for example \cite{ambrosio_metric_1990, KS93, Res97, Haj96, HKST15}. We will review that of Reshetnyak \cite{Res97}.

Let $(X, d)$ be a complete metric space and $M$ a smooth compact $m$-dimen\-sional manifold equipped with a Riemannian metric $g$, possibly with non-empty boundary. Let $ U\subset M$ be an open set and $p>1$. We denote by $L^p( U, X)$ the collection of measurable and essentially separably valued maps $u\colon  U\to X$ such that for some and thus every $x\in X$ the function $u_x\colon U\to \R$, defined by $$u_x(z)= d(x, u(z)),$$ belongs to the classical space $L^p( U)$.

\bd\label{def:Sobolev}
    A map $u\in L^p( U, X)$ belongs to the Sobolev space $W^{1,p}( U, X)$ if for every $x\in X$ the function $u_x$ belongs to $W^{1,p}( U\setminus\partial M)$ and there exists $h\in L^p( U)$ such that for all $x\in X$ we have $|\nabla u_x|_g\leq h$ almost everywhere on $ U$.  
\ed

In the definition above, $\nabla u_x$ is the weak gradient of $u_x$. Further, $|\cdot|_g$ denotes the norm induced by $g$. The Reshetnyak $p$-energy of $u\in W^{1,p}( U, X)$ with regards to $g$ is defined by $$E_+^p(u,g)= \inf_h\int_U h^p(x)\,d\hm_g^m(x),$$ where the infimum is taken over all functions $h$ as in \Cref{def:Sobolev}.

Finally, we recall the definition of the trace of a Sobolev map. Assume $M$ is a smooth surface and let $ U \subset M\setminus\partial M$ be a Lipschitz domain. Then for every $z$ in the boundary $\partial  U$ of $ U$ there exist an open neighborhood $ U_z\subset M$ and a biLip\-schitz map $\psi\colon (0,1)\times [0,1)\to M$ such that $\psi((0,1)\times (0,1)) =  U_z\cap  U$ and $\psi((0,1)\times\{0\}) =  U_z\cap \partial U$. Let $u\in W^{1,2}( U, X)$. For almost every $s\in (0,1)$ the map $t\mapsto u\circ\psi(s,t)$ belongs to $W^{1,2}((0,1),X)$ by a Fubini-type argument and thus has an absolutely continuous representative which we denote by the same expression. The trace of $u$ is defined by $$\trace(u)(\psi(s,0))=\lim_{t\searrow 0} (u\circ\psi)(s,t)$$ for almost every $s\in(0,1)$. It can be shown that the trace is independent of the choice of the map $\psi$ and that it defines an element of $L^2(\partial  U, X)$, see \cite[Section 1.12]{KS93} for further details.

\subsection{Metric differentiability and infinitesimal isotropy}\label{sec:metric-diff}
Let $(X, d)$ be a complete metric space and $M$ a smooth compact surface equipped with a Riemannian metric $g$. Let $U\subset M$ be open and let $p>1$. Sobolev maps $u\in W^{1,p}( U, X)$ have the following approximate metric differentiability property. For almost every $z\in  U$ there exists a unique seminorm $\apmd u_z$ on $T_zM$ such that $$\ap \lim_{v\to 0} \frac{d(u(\exp_z(v)), u(z)) - \apmd u_z(v)}{|v|_g}=0.$$ Here, $\exp_z\colon T_zM\to M$ denotes the exponential map, and $\ap\lim$ is the approximate limit, see \cite{Kir94} and e.g.\ \cite[Theorem 1.15 and Property 2.7]{Kar07}.

Recall that, by John’s theorem \cite{John48}, the unit ball with respect to a norm $||\cdot||$ on $\R^2$ contains a unique ellipse of maximal area, called the John’s ellipse of $||\cdot||$.
\bd\label{def:inf-isotropy}
    A map $u\in W^{1,2}( U, X)$ is infinitesimally isotropic (with respect to $g$) if for almost every $z\in U$ the approximate metric derivative $\apmd u_z$ is either zero or it is a norm and the John's ellipse of $\apmd u_z$ is a round ball (with respect to $g$).
\ed
We say that $X$ has property (ET) if for every $u\in W^{1,2}(U, X)$ the approximate metric derivative $\apmd u_z$ is induced by a possibly degenerate inner product at almost every $z\in U$. Examples of such spaces are Riemannian manifolds with continuous metric tensor, metric spaces of curvature bounded from above or below in the sense of Alexandrov, equiregular sub-Riemannian manifolds, and many more, see \cite[Section 11]{LW-area}.
\bl\label{lem:inf-isotropic-implies-wqc}
    Every monotone and infinitesimally isotropic map $u\in W^{1,2}( U,X)$ satisfies 
    \begin{equation*}
        \MOD(\Gamma)\leq K\cdot \MOD(u\circ\Gamma)
    \end{equation*}
    for every family $\Gamma$ of curves in $U$, where $K=4/\pi$. If $X$ satisfies property (ET), then the statement holds with $K=1$.
\el
\begin{proof}
    Let $z\in U$ be such that $\apmd u_z$ is a norm and denote by $B_z$ the unit ball of $(\R^2,\apmd u_z)$. By John's theorem \cite{John48} and the roundness of the John's ellipse $E_z$, we obtain $\mathcal{L}^2(B_z)\leq K\cdot\mathcal{L}^2(E_z)$ for $K=4/\pi$. This and the roundness of $E_z$ imply that $u$ is infinitesimally $K$-quasiconformal as defined in \cite[Equation (3)]{LW20}. The modulus inequality with constant $K$ now follows from \cite[Proposition 3.5]{LW20}.
    
    If $X$ satisfies property (ET), then we may choose $K=1$ in the arguments above. Indeed, if $\apmd u_z$ is induced by an inner product, then $B_z$ is an ellipse and thus coincides with $E_z$. The roundness of $E_z$ now implies infinitesimal 1-quasiconformality of $u$.
\end{proof}

\section{Monotone Sobolev extension theorems for metric surfaces}\label{sec:Sobolev-monotone-extension}
This section is devoted to proving our main result, \Cref{thm:monotone-filling}. Let $X$ be a metric surface with locally finite Hausdorff 2-measure. We first present the proof under the additional assumption that $X$ is locally geodesic.  Consider a Jordan domain  $\Omega\subset X$  with rectifiable boundary. By scaling the metric on $X$, we may assume that $\partial\Omega$ has length $2\pi$. Let $\alpha\colon S^1\to\partial\Omega$ be a unit speed parametrization and let $Y$ be the metric space obtained by attaching the cylinder $S^1\times[0,2\pi]$ to $\overline{\Omega}$ along $\partial\Omega$ via the identification $(t, 0)\sim\alpha(t)$, and equip $Y$ with the intrinsic metric. It follows from \cite[Lemma 2.1]{LW20} that $Y$ is homeomorphic to $\overline D$ and that $\mathcal H^2(Y)<\infty$. Furthermore, note that the subspace metric of $S^1\times[\pi,2\pi]$ inherited from $Y$ coincides with the Euclidean intrinsic metric.
    
By \cite[Theorem 1.4]{MW21}, the family $\Lambda(\partial Y,Y)$ of Sobolev maps in $W^{1,2}(D, Y)$ whose trace is a monotone parametrization of $\partial Y$ is non-empty. Further, \cite[Theorem 7.6]{LW-area} implies the existence of a Reshetnyak energy minimizer $v\in \Lambda(\partial Y,Y)$. By \cite[Theorem 1.3]{MW21}, the energy minimizer has a continuous representative that extends continuously to the boundary. We denote this representative again by $v\colon\overline D\to Y$. It follows from \cite[Theorem 1.2]{LW20} that $v$ is monotone. By \cite[Lemmas 4.1 and 3.2]{lytchak_energy_2017}, the map $v$ is infinitesimally isotropic. 

We now make use of the local Euclidean structure of $Y$ in a neighborhood of its boundary to prove the following. 
\bl \label{lma:v-smooth-on-the-boundary}
    The map $v|_{S^1}$ is smooth, when considered as a map from $S^1$ to itself. 
\el
\begin{proof}
    Choose $A\subset D$ to be an open annulus with $S^1\subset\partial A$ and $v(A)\subset S^1\times [\pi,2\pi]$. First note that by arguing as in the proof of \cite[Theorem 3.6]{LW20} it holds that $v|_{\overline A}$ is injective. It follows, in particular, that $v|_A$ is a homeomorphism onto its image. As $v|_{A}$ has image in a subspace of $Y$ that satisfies property (ET), we know from \Cref{lem:inf-isotropic-implies-wqc} that
    \begin{equation}\label{ineq:modulus}
        \MOD(\Gamma)\leq\MOD(v\circ\Gamma)
    \end{equation}
    holds for every family $\Gamma$ of curves in $A$. Every homeomorphism between planar Euclidean domains satisfying \eqref{ineq:modulus} also satisfies the reverse inequality of \eqref{ineq:modulus}, see e.g.\ \cite[Section I.3]{Lehto-Virtanen73}. As a consequence, the map $v|_A$ is conformal. It now follows from the Schwarz reflection principle, e.g.\ found in \cite[Theorem 24]{ahlfors_complex_1979}, that $v|_{A}$ is smooth up to the boundary $S^1$.
\end{proof}
We also need the next lemma that provides the existence of a monotone Sobolev homotopy between monotone Sobolev loops in $S^1$.
\bl \label{lma:H-monotone-sobolev}
    Let $f,g\colon S^1\to S^1$ be monotone and homotopic. Then there exists a monotone homotopy $H\colon S^1\times[0,1]\to S^1$ from $f$ to $g$. Moreover, if $f,g\in W^{1,2}(S^1,S^1)$ then we may also assume that $H\in W^{1,2}(S^1\times[0,1],S^1)$.
\el

\begin{proof}
After precomposing $g$ with a rotation of angle $\theta$, we may assume that there is some $z_0\in S^1$ such that $f(z_0)=g(z_0)$. Then, viewing $S^1$ as the interval $[0,1]$ with identified endpoints the maps $f$ and $g$ can be seen as continuous and increasing functions $f, g\colon[0,1]\to[0,1]$ such that $f(0)=  g(0)=0$ and $ f(1)= g(1)=1$. We let $H\colon[0,1]^2\to [0,1]$ be the straight line homotopy between $f$ and $g$, i.e.\ $$H(s,t) = (1-t) f(s)+t g(s).$$

We now show that $H$ is monotone. Note that $H=\pi_1\circ \tilde H$ where $\tilde H\colon [0,1]^2\to [0,1]^2$ is defined by $\tilde H(s,t) = (H(s,t),t)$ and $\pi_1\colon [0,1]^2\to [0,1]$ is the projection map to the first factor. Note that $\pi_1$ and $\tilde H$ are surjective, and that $\pi_1$ is monotone. To see that $\tilde H$ is monotone, observe that for $w,t\in [0,1]$ we have
$$
    \tilde H^{-1}(\{(w,t)\}) = \{(s,t): s\in[0,1], \ (1-t)f(s)+tg(s) = w\},
$$
which is connected since $(1-t)f(\cdot)+tg(\cdot)$ is increasing and hence monotone for every $t\in [0,1]$. It now follows that $H$ is monotone by \cite[Chap. VIII (3.5)]{whyburn_analytic_1942}. By identifying endpoints, and precomposing with a rotation of angle $\theta\cdot t$, the homotopy $H$ descends to a monotone homotopy  $H\colon S^1\times[0,1]\to S^1$ from $f$ to $g$ (now again defined on $S^1)$.

In the case when $f,g\in W^{1,2}(S^1,S^1)$, it follows by a simple computation of the weak partial derivatives that $H\in W^{1,2}(S^1\times[0,1],S^1)$.
\end{proof}

We now return to the proof of \Cref{thm:monotone-filling}. Denote by $\gamma\colon S^1\to\partial\Omega$ the given monotone Sobolev boundary map and let $\alpha$ be the unit speed parametrization of $\partial\Omega$ chosen at the beginning of the current section. Note that $\alpha^{-1}\colon\partial\Omega\to S^1$ preserves lengths of curves and thus $\alpha^{-1}\circ\gamma\in W^{1,2}(S^1,S^1)$. 

We write $v|_{S^1}=(\beta,2\pi)$, where $\beta\colon S^1\to S^1$. We may assume that $\beta$ is homotopic to $\alpha^{-1}\circ \gamma$. By \Cref{lma:H-monotone-sobolev}, there is a monotone homotopy $H\colon S^1\times [0,1]\to S^1$ from $\beta$ to  $\alpha^{-1}\circ \gamma$. Also, note that we may assume that $H$ is Sobolev since by \Cref{lma:v-smooth-on-the-boundary}, $\beta$ is smooth and further $\alpha^{-1}\circ\gamma$ is Sobolev. 
    
Let $\pi\colon Y\to \overline\Omega$ be the Lipschitz map defined such that $\pi|_{\overline\Omega}$ is the identity map on $\overline \Omega$, and $\pi(s,t) \coloneqq \alpha(s)$ for $(s,t)\in S^1\times[0,2\pi]$. The map $u\colon\overline D\to \overline\Omega$ is defined by
\begin{gather*}
    u(z)\coloneqq \begin{cases}
        \pi\circ v(2z),   &|z|\leq 1/2, \\
        \alpha\circ  H(\frac{z}{|z|}, 2|z|-1),  &|z|>1/2.
    \end{cases}
\end{gather*}
It follows that $u\in W^{1,2}(D,\overline\Omega)$ by the Sobolev gluing theorem, see \cite[Theorem 1.12.3]{KS93}. Further, note that if $z\in S^1$, then $u(z) = \gamma(z)$.

It remains to show that $u$ is monotone. Monotonicity of $v$ directly implies that $u^{-1}(\{x\})$ is connected for $x\in\Omega$. Now, consider a value from the boundary $x\in\partial\Omega$ and let $x'\coloneqq\alpha^{-1}(x)$. Then $u^{-1}(\{x\}) = v^{-1}(T)\cup H^{-1}(\{x'\})$, where $T$ is the interval $\{x'\}\times[0,2\pi]\subset S^1\times[0,2\pi]$. Since $\overline D$ is compact and $v$ monotone, it follows from \cite[Chap. VIII (2.2)]{whyburn_analytic_1942} that $v^{-1}(T)$ is connected. Furthermore, since $H$ is monotone, $H^{-1}(\{x'\})$ is connected. The two preimages in the union contain a common point, and thus, $v^{-1}(T)\cup H^{-1}(\{x'\})$ is connected. We conclude that $u$ is monotone. This completes the proof of \Cref{thm:monotone-filling}  in the case where $X$ is locally geodesic. By exploiting the more recent work \cite{NR22} instead of \cite{MW21}, we may prove the theorem in full generality. This is explained in more detail in the next remark. 

\br\label{rmk:NR22}
    Assume $X$ is a metric surface of locally finite Hausdorff 2-measure that is not necessarily locally geodesic. We may construct the space $Y$ as above and equip it with the quotient metric $d_q$; for a definition we refer to \cite[Definition 3.1.12]{BBI01}. Note that the space $(Y,d_q)$ is homeomorphic to $\overline{D}$ and the natural quotient map from $(S^1\times[0,2\pi])\sqcup\overline\Omega$ to $Y$ is a local isometry on $\Omega$ and on $S^1\times (0,2\pi]$. In particular, $\hm^2(Y)<\infty$ and we may apply \cite[Theorem~1.2]{NR22} to get the existence of a weakly $(4/\pi)$-quasiconformal parametrization $v\colon\overline{D}\to Y$. By definition, $v$ is monotone and surjective. Moreover, it follows from the proof that $v$ is infinitesimally isotropic, see \cite[Proof of Theorem 6.1]{NR22}. These are all the properties of $v$ needed for continuing the proof as above and thus, arriving at the conclusion of \Cref{thm:monotone-filling}. 
\er

We end this section by noting that,  in the setting of the theorem, there is no assumption we can make on the boundary map that would guarantee the existence of a homeomorphic Sobolev extension. This is illustrated in the following example.
\begin{example}\label{ex:collapsed-disc}  
    Let $X$ be the space obtained from $\R^2$ by collapsing a disc of radius $s_0>0$ to a point $p_0$. We endow $X$ with the quotient metric, see \cite[Definition 3.1.12]{BBI01}.  Then, the class of Sobolev homeomorphisms from $\overline{D}$ to $\overline{\Omega}$ is empty, where $\Omega\subset X$ is some Jordan domain containing $p_0$. Indeed, assume for contradiction that there exists a Sobolev homeomorphism $u\colon\overline D\to \overline\Omega$, and let $z_0=u^{-1}(p_0)$. Then, for any $r\in(0,1-|z_0|)$ we have that  $\ell(u|_{\partial B(z_0,r)})\geq 2\pi s_0$, since $u|_{\partial B(z_0,r)}$ is a loop around $p_0$. However, by the Courant-Lebesgue lemma, see e.g.\ \cite[Lemma 7.3]{LW-area}, we can choose $r\in(0,1-|z_0|)$ to make $\ell(u|_{\partial B(z_0,r)})$  arbitrarily small, which stands in contradiction to the previous sentence. 
\end{example}

\subsection{Monotone Sobolev extensions to surfaces of higher topology}
Using \Cref{thm:monotone-filling}, we can easily construct monotone Sobolev extensions defined on compact smooth surfaces.

\bt \label{thm:monotone-sobolev-extension-compact-surface}
    Let $X$ be a metric surface with locally finite Hausdorff 2-measure. Let $\overline \Omega$ be a subset of $X$ that is homeomorphic to a compact smooth surface $\Sigma$ with non-empty boundary such that $\partial \Omega$ has finite length. Then every monotone map in $W^{1,2}(\partial \Sigma,\partial \Omega)$ with a continuous extension from $\Sigma$ to $\overline \Omega$ also extends to a monotone map in $W^{1,2}(\Sigma,\overline \Omega)$. 
\et

\begin{proof}
     The case of $\Sigma$ being a disc readily follows from \Cref{thm:monotone-filling}. Thus, we assume that $\Sigma$ is not of disc-type. For simplicity, we furthermore assume that $\Sigma$ is orientable, the non-orientable case follows by similar arguments.
     Let $\varrho\in W^{1,2}(\partial\Sigma,\partial\Omega)$ be the monotone boundary map. The fact that $\varrho$ extends continuously to $\Sigma$ implies, by \cite{youngsExtensionHomeomorphismDefined1948},  that there also exists a homeomorphism $\varphi\colon\Sigma\to\Omega$  satisfying $d(\varrho(x),\varphi(x))<\varepsilon$ for a fixed small enough $\varepsilon>0$ and every $x\in\partial\Sigma$.

     Choose a collection of simple closed curves $\gamma_i$ decomposing $\Sigma$ into smooth Y-pieces and cylinders $\Sigma_k$. Here, a surface of genus $0$ is called Y-piece if it has three boundary components. By adding smooth curves between boundary components, we conclude that every $\Sigma_k$ decomposes into two sets $\Sigma_{k,1}$ and $\Sigma_{k,2}$ that are biLipschitz equivalent to $\overline{D}$. 

     As in the proof of \cite[Proposition 3.5]{Raj17}, we make use of planar topological arguments and the coarea inequality for Lipschitz functions, see e.g.\ \cite[Theorem~2.10.25]{Fed69}, to obtain the following:
     for every $\gamma_i$ there exists a rectifiable Jordan curve $\eta_i$ homotopic to $\varphi|_{\gamma_i}$ and with image in an arbitrarily small neighborhood of $\varphi(\gamma_i)$. This induces a decomposition of $\overline{\Omega}$ into Y-pieces and cylinders $\Omega_k$ with rectifiable boundary. Here, we choose the indexing such that $\Omega_k$ is contained in the small neighborhood of $\varphi(\Sigma_k)$. Moreover, for every component $\gamma'$ of $\partial\Sigma_{k,1}\cap\partial\Sigma_{k,2}$ there exists a rectifiable Jordan arc $\eta'$ contained in an arbitrarily small neighborhood of $\varphi(\gamma')$ connecting two boundary components of $\Omega_k$. In particular, every $\Omega_k$ decomposes into two sets $\Omega_{k,1}$ and $\Omega_{k,2}$ that have rectifiable boundary and are homeomorphic to $\overline{D}$. Again, we choose the indexing such that $\Omega_{k,l}$ is contained in the small neighborhood of $\varphi(\Sigma_{k,l})$. 

     Denote by $\Sigma^1\subset\Sigma$ and $\Omega^1\subset\overline{\Omega}$ the union of all $\partial\Sigma_{k,l}$ and all $\partial\Omega_{k,l}$, respectively. 
     By a slight adjustment of $\Sigma^1$ close to the boundary $\partial\Sigma$, we may assume that  $$\varrho(\overline{(\Sigma^1\setminus\partial\Sigma)}\cap\partial\Sigma) = \overline{(\Omega^1\setminus\partial\Omega)}\cap\partial\Omega.$$ 
     The above construction allows us to extend $\varrho$ to a monotone map $\chi\colon\Sigma^1\to\Omega^1$ that is Lipschitz on $\Sigma^1\setminus\partial \Sigma$, sends each $\partial\Sigma_{k,l}$ to $\partial\Omega_{k,l}$, and is homotopic to $\varphi|_{\Sigma^1}$ via a homotopy that sends $\partial\Sigma$ to $\partial\Omega$.

     We now define an extension $u\colon \Sigma\to \overline\Omega$ of $\varrho$ piecewise on every $\Sigma_{k,l}\subset \overline\Omega$. By  \Cref{thm:monotone-filling} and the fact that $\Sigma_{k,l}$ is biLipschitz equivalent to $\overline{D}$, we may define $u|_{\Sigma_{k,l}}$ to be a monotone Sobolev extension of $\chi|_{\partial\Sigma_{k,l}}$ such that $u(\Sigma_{k,l})=\Omega_{k,l}$. Then, $u$ is a Sobolev map by the Sobolev gluing theorem, see \cite[Theorem 1.12.3]{KS93}. Further, $u$ is clearly surjective. Now, we show that $u$ is monotone. It follows directly from the monotonicity of the restrictions $u|_{\Sigma_{k,l}}$ that fibers of points in $\overline\Omega\setminus \Omega^1$ under $u$ are connected. So, let $x\in \Omega^1$ and observe that the fiber of $x$ can be written as a finite union of connected sets, that is
    $$
    u^{-1}(\{x\})=\bigcup_{k,l}(u|_{\Sigma_{k,l}})^{-1}(\{x\}).
    $$
    Note that $(u|_{\Sigma_{k,l}})^{-1}(\{x\})$ is non-empty for some $i$ only if it contains $g^{-1}(\{x\})$ as a subset, so the non-empty sets in the above union have a common point. Therefore $u^{-1}(\{x\})$ is connected. 
\end{proof}

\br
    The condition in \Cref{thm:monotone-sobolev-extension-compact-surface} that the monotone boundary map admits a continuous extension is not needed for non-orientable surfaces, by \cite{youngsExtensionHomeomorphismDefined1948}.
\er

By adjusting the above arguments to include the cases of closed surfaces, the first statement of \Cref{thm:uniformization} follows as a corollary. 

\bc \label{lma:Lambda-non-empty}
    Let $X$ be a metric surface that is homeomorphic to a compact smooth orientable surface $M$, has rectifiable boundary (if non-empty) and is of finite Hausdorff 2-measure. Then the family $\Lambda$ of surjective and monotone maps in $W^{1,2}(M,X)$ is non-empty. 
\ec

\section{Uniformization of metric surfaces via energy minimization}\label{sec:uniformization}

The goal of this section is to prove our result, \Cref{thm:uniformization}, on the uniformization of metric surfaces via energy minimizers. Recall that in \Cref{thm:uniformization} the metric space $X$ is locally geodesic, homeomorphic to a compact orientable smooth surface $M$, has rectifiable boundary (if non-empty) and is of finite Hausdorff 2-measure.

The monotone Sobolev extension results proved in \Cref{sec:Sobolev-monotone-extension} play an important role in our approach. Indeed, by \Cref{lma:Lambda-non-empty}, the family $\Lambda$ of surjective and monotone Sobolev maps from $M$ to $X$ is non-empty. Thus, we may choose a sequence of pairs $(v_n,h_n)$ such that
$$
E_+^2(v_n,h_n)\to e(\Lambda)\coloneqq\inf\{E_+^2(v,h): v\in \Lambda, \ h \text{ Riemannian metric on $M$}\},
$$
where $v_n\in \Lambda$ and $h_n$ is a Riemannian metric on $M$. By the classical uniformization theorem, we may assume that the metrics $h_n$ are of constant sectional curvature $-1,0$ or $1$ and such that each component of $\partial M$ is geodesic (if non-empty). If $M$ admits metrics of curvature $0$, then we normalize $h_n$ such that $\text{Area}(M,h_n)=1$.

Assume first that $M$ is not simply connected. In this situation, \cite[Proposition 8.4]{fitziAreaMinimizingSurfaces2021}, whose proof can be adapted to include the easier cases of $M$ being a torus or cylinder, implies that there is a uniform lower bound on relative systoles of $(M,h_n)$. By Mumford's compactness theorem, see e.g.\ \cite[Theorem 3.3]{fitziAreaMinimizingSurfaces2021}, it follows that there exist diffeomorphisms $\varphi_n\colon M\to M$ such that, up to taking a subsequence, $\varphi_n^*h_n$ converges smoothly to some Riemannian metric $g$ of constant sectional curvature. Let $u_n\coloneqq v_n\circ \varphi_n\in \Lambda$ and note that $E^2_+(u_n,g)$ converges to $e(\Lambda)$. This is because for every $\varepsilon>0$ the map $\varphi_n\colon(M,g)\to (M,h_n)$ is $(1+\varepsilon)$-biLipschitz for large $n$.

Applying \cite[Proposition 4.1]{fitzi_canonical_2022}  implies that $(u_n)$ is equicontinuous. Hence, we know that a subsequence converges uniformly to a continuous and surjective map $u\colon M\to X$. This map is Sobolev as a consequence of Rellich-Kondrachov compactness, see \cite[Theorem 1.13]{KS93}. Further, since $u$ is a uniform limit of monotone maps valued in a locally connected space, it follows that also $u$ is monotone by \cite[Chap. IX, Corollary 3.11]{whyburn_analytic_1942}. Hence, it is immediate that $u\in\Lambda$. By the lower-semicontinuity of energy, we have $E^2_+(u,g)=e(\Lambda)$. This fact implies that $u$ is infinitesimally isotropic with respect to $g$ by Morrey's $\varepsilon$-conformality Lemma, see \cite[Corollary 1.3]{fitziMorreysvarepsilonconformality2020}. Lastly, \Cref{lem:inf-isotropic-implies-wqc} shows that $u$ is weakly $4/\pi$-quasiconformal. This concludes the proof in the case that $M$ is not simply connected.

Now, we will describe how to adapt the above proof to the case when $M$ is simply connected, or equivalently, when $M$ is homeomorphic to either the closed disc or the sphere. In this situation, the Mumford's compactness theorem is not needed. This is because all Riemannian metrics on $M$ are conformally equivalent, meaning that we may assume that $h_n$ are all equal to the same Riemannian metric $g$ of constant curvature. 

Up to precompositions of Möbius transformations, we may assume that $(v_n)$ satisfies a three-point condition. Then, equicontinuity of the sequence follows by arguing as in the proof of \cite[Lemma 3.1.2]{Jost_two-dim} for the sphere and \cite[Theorem 4.3.1]{dierkes_minimal_2010} for the disc, using an appropriate variant of the Courant-Lebesgue Lemma, see \cite[Lemma 7.3]{LW-area}. The rest of the argument from the non-simply connected case can be applied verbatim from this point on. This completes the proof of \Cref{thm:uniformization}.

\bibliographystyle{plain}
\bibliographystyle{plain}
    
\end{document}